\documentclass{amsart}
\usepackage{amssymb}
\usepackage{amsfonts}

\newtheorem{theorem}{Theorem}
\theoremstyle{plain}

\newtheorem{corollary}{Corollary}

\newtheorem{proposition}{Proposition}
\newtheorem{remark}{Remark}

\numberwithin{equation}{section}
\input{tcilatex}

\begin{document}
\title[Schwarz's Inequality]{Refinements of Some Reverses of Schwarz's
Inequality in $2-$Inner Product Spaces and Applications for Integrals}
\author{P. Cerone}
\address{School of Computer Science and Mathematics\\
Victoria University of Technology\\
PO Box 14428, MCMC 8001\\
Victoria, Australia.}
\email{pietro.cerone@vu.edu.au}
\urladdr{http://rgmia.vu.edu.au/cerone}
\author{Y.J. Cho$^{\bigstar }$}
\address{Department of Mathematics Education, The Research Institute of
Natural Sciences, Gyeongsang National University\\
Chinju 660-701, Korea}
\email{yjcho@nongae.gsnu.ac.kr}
\author{S.S. Dragomir}
\address{School of Computer Science and Mathematics\\
Victoria University of Technology\\
PO Box 14428, MCMC 8001\\
Victoria, Australia.}
\email{sever.dragomir@vu.edu.au}
\urladdr{http://rgmia.vu.edu.au/SSDragomirWeb.html}
\author{S.S. Kim$^{\blacklozenge }$}
\address{Department of Mathematics, Dongeui University\\
Pusan 614-714, Korea}
\email{sskim@dongeui.ac.kr}
\date{August 28, 2003}
\keywords{$2-$Inner product spaces, Schwarz's inequality, Determinantal
integral inequalities\\
\ \ \ \ \ \ \ $\bigstar ,\blacklozenge $ Corresponding authors}

\begin{abstract}
Refinements of some recent reverse inequalities for the celebrated
Cauchy-Bunyakovsky-Schwarz inequality in $2-$inner product spaces are given.
Using this framework, applications for determinantal integral inequalities
are also provided.
\end{abstract}

\maketitle

\section{Introduction}

The concepts of $2-$inner products and $2-$inner product spaces have been
intensively studied by many authors in the last three decades.

A systematic presentation of the recent results related to the theory of $2-$%
inner product spaces as well as an extensive list of the related references
can be found in the book \cite{CLKM}. We recall here the basic definitions
and the elementary properties of $2-$inner product spaces that will be used
in the sequel (see also \cite{CMP1}).

Let $X$ be a linear space of dimension greater than $1$ over the number
field $\mathbb{K}$, when $\mathbb{K}=\mathbb{R}$ or $\mathbb{K}=\mathbb{C}$.
Suppose that $\left( \cdot ,\cdot |\cdot \right) $ is a $\mathbb{K}$-valued
function defined on $X\times X\times X$ satisfying the following conditions:

\begin{enumerate}
\item[$\left( 2I_{1}\right) $] $\left( x,x|z\right) \geq 0$ and $\left(
x,x|z\right) =0$ if and only if $x$ and $z$ are linearly dependent,

\item[$\left( 2I_{2}\right) $] $\left( x,x|z\right) =\left( z,z|x\right) ,$

\item[$\left( 2I_{3}\right) $] $\left( y,x|z\right) =\overline{\left(
x,y|z\right) },$

\item[$\left( 2I_{4}\right) $] $\left( \alpha x,y|z\right) =\alpha \left(
x,y|z\right) $ for any scalar $\alpha \in \mathbb{K}$,

\item[$\left( 2I_{5}\right) $] $\left( x+x^{\prime },y|z\right) =\left(
x,y|z\right) +\left( x^{\prime },y|z\right) ,$
\end{enumerate}

where $x,x^{\prime },y,z\in X.$ The functional $\left( \cdot ,\cdot |\cdot
\right) $ is called a $2-$\textit{inner product }on $X$ and $\left( X,\left(
\cdot ,\cdot |\cdot \right) \right) $ is called a $2-$\textit{inner product
space} (or \textit{2-pre-Hilbert space}) \cite{CLKM}.

Some basic properties of the $2-$inner product spaces can be immediately
obtained as follows:

\begin{enumerate}
\item If $\mathbb{K}=\mathbb{R}$, then $\left( 2I_{3}\right) $ reduces to%
\begin{equation*}
\left( y,x|z\right) =\left( x,y|z\right) .
\end{equation*}

\item From $\left( 2I_{3}\right) $ and $\left( 2I_{4}\right) ,$ we have%
\begin{equation*}
\left( 0,y|z\right) =\left( x,0|z\right) =0
\end{equation*}%
and also%
\begin{equation}
\left( x,\alpha y|z\right) =\bar{\alpha}\left( x,y|z\right) .  \label{1.1}
\end{equation}

\item Using $\left( 2I_{3}\right) -\left( 2I_{5}\right) ,$ we have%
\begin{eqnarray*}
\left( z,z|x\pm y\right) &=&\left( x\pm y,x\pm y|z\right) \\
&=&\left( x,x|z\right) +\left( y,y|z\right) \pm 2\func{Re}\left( x,y|z\right)
\end{eqnarray*}%
and%
\begin{equation}
\func{Re}\left( x,y|z\right) =\frac{1}{4}\left[ \left( z,z|x+y\right)
-\left( z,z|x-y\right) \right] .  \label{1.2}
\end{equation}%
In the real case $\mathbb{K}=\mathbb{R}$, (\ref{1.2}) reduces to%
\begin{equation}
\left( x,y|z\right) =\frac{1}{4}\left[ \left( z,z|x+y\right) -\left(
z,z|x-y\right) \right] ,  \label{1.3}
\end{equation}%
and using this formula, it is easy to see, for any $\alpha \in \mathbb{R}$,
that%
\begin{equation}
\left( x,y|\alpha z\right) =\alpha ^{2}\left( x,y|z\right) .  \label{1.4}
\end{equation}%
In the complex case, $\mathbb{K}=\mathbb{C}$, using (\ref{1.1}) and (\ref%
{1.2}), we have%
\begin{equation*}
\func{Im}\left( x,y|z\right) =\func{Re}\left[ -i\left( x,y|z\right) \right] =%
\frac{1}{4}\left[ \left( z,z|x+iy\right) -\left( z,z|x-iy\right) \right] ,
\end{equation*}%
which, in combination with (\ref{1.2}), yields%
\begin{equation}
\left( x,y|z\right) =\frac{1}{4}\left[ \left( z,z|x+y\right) -\left(
z,z|x-y\right) \right] +\frac{i}{4}\left[ \left( z,z|x+iy\right) -\left(
z,z|x-iy\right) \right] .  \label{1.5}
\end{equation}%
Using (\ref{1.5}) and (\ref{1.1}), we have, for any $\alpha \in \mathbb{C}$,
that%
\begin{equation}
\left( x,y|\alpha z\right) =\left\vert \alpha \right\vert ^{2}\left(
x,y|z\right) .  \label{1.6}
\end{equation}%
However, for $\alpha \in \mathbb{R}$, (\ref{1.6}) reduces to (\ref{1.4}).
Also, from (\ref{1.6}) it follows that%
\begin{equation*}
\left( x,y|0\right) =0.
\end{equation*}

\item For any three given vectors $x,y,z\in X,$ consider the vector $%
u=\left( y,y|z\right) x-\left( x,y|z\right) y.$ By $\left( 2I_{1}\right) ,$
we know that $\left( u,u|z\right) \geq 0$ with the equality if and only if $%
u $ and $z$ are linearly dependent. It is obvious that the inequality $%
\left( u,u|z\right) \geq 0$ can be rewritten as%
\begin{equation}
\left( y,y|z\right) \left[ \left( x,x|z\right) \left( y,y|z\right)
-\left\vert \left( x,y|z\right) \right\vert ^{2}\right] \geq 0.  \label{1.7}
\end{equation}%
For $x=z,$ (\ref{1.7}) becomes%
\begin{equation*}
-\left( y,y|z\right) \left\vert \left( z,y|z\right) \right\vert ^{2}\geq 0
\end{equation*}%
which implies that%
\begin{equation}
\left( z,y|z\right) =\left( y,z|z\right) =0,  \label{1.8}
\end{equation}%
provided \thinspace $y$ and $z$ are linearly independent. Obviously, when $y$
and $z$ are linearly dependent, (\ref{1.8}) also holds.

Now, if $y$ and $z$ are linearly independent, then $\left( y,y|z\right) >0,$
and from (\ref{1.7}), it follows the Cauchy-Bunyakovsky-Schwarz inequality ($%
CBS-$inequality for short) for $2-$inner products:%
\begin{equation}
\left\vert \left( x,y|z\right) \right\vert ^{2}\leq \left( x,x|z\right)
\left( y,y|z\right) .  \label{1.9}
\end{equation}

Utilizing (\ref{1.8}), it is easy to see that (\ref{1.9}) is trivially
fulfilled when $y$ and $z$ are linearly dependent. Therefore, the inequality
(\ref{1.9}) holds for any three vectors $x,y,z\in X$ and is strict unless
the vectors 
\begin{equation*}
u=\left( y,y|z\right) x-\left( x,y|z\right) y\text{ \ and \ }z
\end{equation*}%
are linearly dependent. In fact, \textit{we have the equality in (\ref{1.9})
if and only if the three vectors }$x,y$ \textit{and }$z$ \textit{are
linearly dependent \cite{CMP1}.}
\end{enumerate}

In any given $2-$inner product space $\left( X,\left( \cdot ,\cdot |\cdot
\right) \right) ,$ we can define a function $\left\Vert \cdot |\cdot
\right\Vert $ on $X\times X$ by%
\begin{equation}
\left\Vert x|z\right\Vert =\sqrt{\left( x,x|z\right) }  \label{1.10}
\end{equation}%
for all $x,z\in X.$ It is easy to see that, this function satisfies the
following conditions

\begin{enumerate}
\item[$\left( 2N_{1}\right) $] $\left\Vert x|z\right\Vert \geq 0$ and $%
\left\Vert x|z\right\Vert =0$ if and only if $x$ and $z$ are linearly
dependent,

\item[$\left( 2N_{2}\right) $] $\left\Vert z|x\right\Vert =\left\Vert
x|z\right\Vert ,$

\item[$\left( 2N_{3}\right) $] $\left\Vert \alpha x|z\right\Vert =\left\vert
\alpha \right\vert \left\Vert x|z\right\Vert $ for any scalar $\alpha \in 
\mathbb{K}$,

\item[$\left( 2N_{4}\right) $] $\left\Vert x+x^{\prime }|z\right\Vert \leq
\left\Vert x|z\right\Vert +\left\Vert x^{\prime }|z\right\Vert .$
\end{enumerate}

Any function $\left\Vert \cdot |\cdot \right\Vert $ defined on $X\times X$
and satisfying the conditions $\left( 2N_{1}\right) -\left( 2N_{4}\right) $
is called a 2-\textit{norm} on $X$ and $\left( X,\left\Vert \cdot |\cdot
\right\Vert \right) $ is called a \textit{linear 2-normed space \cite{FC}.}

In terms of 2-norms, the $\left( CBS\right) -$inequality (\ref{1.9}) can be
written as 
\begin{equation}
\left\vert \left( x,y|z\right) \right\vert ^{2}\leq \left\Vert
x|z\right\Vert ^{2}\left\Vert y|z\right\Vert ^{2}.  \label{1.11}
\end{equation}%
The equality in (\ref{1.11}) holds if and only if $x,y$ and $z$ are linearly
dependent.

For recent inequalities in 2-inner products, see the recent works \cite{BCMP}
- \cite{CSL} and the references therein.

In \cite{DCK}, the authors pointed out the following reverses of the $\left(
CBS\right) -$inequality in 2-inner product spaces.

Assume that $x,y,z\in X$ and $a,A\in \mathbb{K}$ are such that either 
\begin{equation}
\func{Re}\left( Ay-x,x-ay|z\right) \geq 0  \label{1.12}
\end{equation}%
or, equivalently%
\begin{equation}
\left\Vert x-\frac{a+A}{2},y|z\right\Vert \leq \frac{1}{2}\left\vert
A-a\right\vert \left\Vert y|z\right\Vert  \label{1.13}
\end{equation}%
hold. Then one has the inequality \cite{DCK}%
\begin{equation}
0\leq \left\Vert x|z\right\Vert ^{2}\left\Vert y|z\right\Vert
^{2}-\left\vert \left( x,y|z\right) \right\vert ^{2}\leq \frac{1}{4}%
\left\vert A-a\right\vert ^{2}\left\Vert y|z\right\Vert ^{4}.  \label{1.14}
\end{equation}%
The constant $\frac{1}{4}$ is sharp in\ (\ref{1.14}) in the sense that it
cannot be replaced by a smaller constant.

With the same assumptions for $x,y,z,a$ and $A$ and, if moreover $\func{Re}%
\left( \bar{a}A\right) >0$, then \cite{DCK}%
\begin{eqnarray}
\left\Vert x|z\right\Vert \left\Vert y|z\right\Vert  &\leq &\frac{1}{2}\cdot 
\frac{\func{Re}\left[ \left( \bar{A}+\bar{a}\right) \left( x,y|z\right) %
\right] }{\func{Re}\left[ \left( \bar{a}A\right) \right] ^{\frac{1}{2}}}
\label{1.15} \\
&\leq &\frac{1}{2}\cdot \frac{\left\vert A+a\right\vert }{\func{Re}\left[
\left( \bar{a}A\right) \right] ^{\frac{1}{2}}}\left\vert \left( x,y|z\right)
\right\vert .  \notag
\end{eqnarray}%
Here the constant $\frac{1}{2}$ is best possible in both inequalities.

As a consequence of (\ref{1.15}) we may get the following additive reverse
of the $\left( CBS\right) -$inequality as well \cite{DCK}%
\begin{equation}
0\leq \left\Vert x|z\right\Vert ^{2}\left\Vert y|z\right\Vert
^{2}-\left\vert \left( x,y|z\right) \right\vert ^{2}\leq \frac{1}{4}\cdot 
\frac{\left\vert A-a\right\vert ^{2}}{\func{Re}\left( \bar{a}A\right) }%
\left\vert \left( x,y|z\right) \right\vert ^{2}.  \label{1.16}
\end{equation}%
The constant $\frac{1}{4}$ in (\ref{1.16}) is best possible in the above
sense.

\section{Refinements of a Reverse $\left( CBS\right) -$Inequality}

The following reverse of the $\left( CBS\right) -$inequality holds.

\begin{theorem}
\label{t2.1}Let $\left( X,\left( \cdot ,\cdot |\cdot \right) \right) $ be a
2-inner product space on $\mathbb{K}$, $x,y,z\in X$ and $a,A\in \mathbb{K}$.
If%
\begin{equation}
\func{Re}\left( Ay-x,x-ay|z\right) \geq 0,  \label{2.1}
\end{equation}%
or, equivalently,%
\begin{equation}
\left\Vert x-\frac{a+A}{2}y|z\right\Vert \leq \frac{1}{2}\left\vert
A-a\right\vert \left\Vert y|z\right\Vert ,  \label{2.2}
\end{equation}%
holds, then one has the inequality%
\begin{eqnarray}
0 &\leq &\left\Vert x|z\right\Vert ^{2}\left\Vert y|z\right\Vert
^{2}-\left\vert \left( x,y|z\right) \right\vert ^{2}  \label{2.3} \\
&\leq &\frac{1}{4}\left\vert A-a\right\vert ^{2}\left\Vert y|z\right\Vert
^{4}-\left\vert \frac{a+A}{2}\left\Vert y|z\right\Vert ^{2}-\left(
x,y|z\right) \right\vert ^{2}  \notag \\
&&\left( \leq \frac{1}{4}\left\vert A-a\right\vert ^{2}\left\Vert
y|z\right\Vert ^{4}\right) .  \notag
\end{eqnarray}%
The constant $\frac{1}{4}$ is sharp in (\ref{2.3}) in the sense that it
cannot be replaced by a smaller constant.
\end{theorem}

\begin{proof}
Observe , for $x,u,U\in X,$ that we have%
\begin{align*}
\frac{1}{4}\left\Vert U-u|z\right\Vert ^{2}-\left\Vert \left. x-\frac{u+U}{2}%
\right\vert z\right\Vert ^{2}& =\func{Re}\left( U-u,x-u|z\right) \\
& =\func{Re}\left[ \left( x,u|z\right) +\left( U,x|z\right) \right] -\func{Re%
}\left( U,u|z\right) -\left\Vert x,z\right\Vert ^{2}.
\end{align*}%
Therefore%
\begin{equation*}
\func{Re}\left( U-u,x-u|z\right) \geq 0,
\end{equation*}%
if and only if%
\begin{equation*}
\left\Vert \left. x-\frac{u+U}{2}\right\vert z\right\Vert \leq \frac{1}{2}%
\left\Vert U-u|z\right\Vert .
\end{equation*}%
If we choose above $U=Ay$ and $u=ay,$ we deduce that the conditions (\ref%
{2.1}) and (\ref{2.3}) are equivalent.

Now, if we consider $x,y,z\in X$ and $\lambda \in \mathbb{K}$, then we may
state that%
\begin{equation}
\left\Vert x-\lambda y|z\right\Vert ^{2}=\left\Vert x|z\right\Vert ^{2}-2%
\func{Re}\left[ \lambda \left( x,y|z\right) \right] +\left\vert \lambda
\right\vert ^{2}\left\Vert y|z\right\Vert ^{2}  \label{2.4}
\end{equation}%
and 
\begin{equation}
\left\vert \lambda \left\Vert y|z\right\Vert ^{2}-\left( x,y|z\right)
\right\vert ^{2}=\left\vert \lambda \right\vert ^{2}\left\Vert
y|z\right\Vert ^{2}-2\left\Vert y|z\right\Vert ^{2}\func{Re}\left[ \lambda
\left( x,y|z\right) \right] +\left\vert \left( x,y|z\right) \right\vert .
\label{2.5}
\end{equation}%
If we multiply (\ref{2.4}) by $\left\Vert x|z\right\Vert ^{2}\geq 0$ and
then subtract equation (\ref{2.5}), we deduce the following equality, that
is of interest in itself,%
\begin{equation}
\left\Vert x|z\right\Vert ^{2}\left\Vert y|z\right\Vert ^{2}-\left\vert
\left( x,y|z\right) \right\vert ^{2}=\left\Vert x-\lambda y|z\right\Vert
^{2}\left\Vert y|z\right\Vert ^{2}-\left\vert \lambda \left\Vert
y|z\right\Vert ^{2}-\left( x,y|z\right) \right\vert ^{2}.  \label{2.6}
\end{equation}%
If we now use (\ref{2.6}) for $\lambda =\frac{a+A}{2}$ and take into account
(\ref{2.2}), then we deduce the desired inequality (\ref{2.3}).

To prove the sharpness of the constant $\frac{1}{4}$ in the second
inequality in (\ref{2.3}), assume that, this inequality holds with a
constant $C>0.$ That is,%
\begin{equation}
\left\Vert x|z\right\Vert ^{2}\left\Vert y|z\right\Vert ^{2}-\left\vert
\left( x,y|z\right) \right\vert ^{2}\leq C\left\vert A-a\right\vert
^{2}\left\Vert y|z\right\Vert ^{4}-\left\vert \frac{a+A}{2}\left\Vert
y|z\right\Vert ^{2}-\left( x,y|z\right) \right\vert ^{2},  \label{2.7}
\end{equation}%
where $x,y,z,a$ and $A$ satisfy the hypothesis of the theorem.

Consider $y,z\in X$ with $\left\Vert y|z\right\Vert =1,$ $a\neq A,$ $a,A\in 
\mathbb{K}$ and $m\in X$ with $\left\Vert m|z\right\Vert =1$ and $\left(
y,m|z\right) =0.$ Define the vector%
\begin{equation*}
x:=\frac{a+A}{2}y+\frac{A-a}{2}m.
\end{equation*}%
Then a simple calculation shows that%
\begin{equation*}
\left( Ay-x,x-ay|z\right) =\frac{\left\vert A-a\right\vert ^{2}}{4}\left(
y-m,y+m|z\right) =0,
\end{equation*}%
and thus the condition (\ref{2.1}) is fulfilled.

Observe also that%
\begin{equation*}
\left\Vert x|z\right\Vert ^{2}=\left\Vert \left. \frac{a+A}{2}y+\frac{A-a}{2}%
m\right\vert z\right\Vert =\left\vert \frac{A+a}{2}\right\vert
^{2}+\left\vert \frac{A-a}{2}\right\vert ^{2},
\end{equation*}%
and%
\begin{equation*}
\left( x,y|z\right) =\left( \frac{a+A}{2}y+\frac{A-a}{2}m,y|z\right) =\frac{%
a+A}{2}.
\end{equation*}%
Consequently, by (\ref{2.7}), we deduce%
\begin{equation*}
\frac{\left( A-a\right) ^{2}}{4}\leq C\left\vert A-a\right\vert ^{2},
\end{equation*}%
giving $C\geq \frac{1}{4},$ and the theorem is proved.
\end{proof}

Another reverse for the (CBS)-inequality is incorporated in the following
theorem.

\begin{theorem}
\label{t2.2}With the assumptions of Theorem \ref{t2.1}, one has the
inequality%
\begin{eqnarray}
0 &\leq &\left\Vert x|z\right\Vert ^{2}\left\Vert y|z\right\Vert
^{2}-\left\vert \left( x,y|z\right) \right\vert ^{2}  \label{2.8} \\
&\leq &\frac{1}{4}\left\vert A-a\right\vert ^{2}\left\Vert y|z\right\Vert
^{4}-\func{Re}\left( Ay-x,x-ay|z\right) \left\Vert y|z\right\Vert ^{2} 
\notag \\
&&\left( \leq \frac{1}{4}\left\vert A-a\right\vert ^{2}\left\Vert
y|z\right\Vert ^{4}\right) .  \notag
\end{eqnarray}%
The constant $\frac{1}{4}$ is sharp in (\ref{2.8}).
\end{theorem}

\begin{proof}
We use the following identity that has been obtained in \cite{DCK} and can
be proved by direct computation%
\begin{multline}
\left\Vert x|z\right\Vert ^{2}\left\Vert y|z\right\Vert ^{2}-\left\vert
\left( x,y|z\right) \right\vert ^{2}  \label{2.9} \\
=\func{Re}\left[ \left( A\left\Vert y|z\right\Vert ^{2}-\left( x,y|z\right)
\right) \left( \overline{\left( x,y|z\right) }-\bar{a}\left\Vert
y|z\right\Vert ^{2}\right) \right] \\
-\left\Vert y|z\right\Vert ^{2}\func{Re}\left( Ay-x,x-ay|z\right) .
\end{multline}%
By the elementary inequality%
\begin{equation*}
\func{Re}\left( \alpha \bar{\beta}\right) \leq \frac{1}{4}\left\vert \alpha
+\beta \right\vert ^{2},\ \ \ \ \alpha ,\beta \in \mathbb{K}
\end{equation*}%
applied for 
\begin{equation*}
\alpha :=A\left\Vert y|z\right\Vert ^{2}-\left( x,y|z\right) \text{ \ and \ }%
\beta =\left( x,y|z\right) -a\left\Vert y|z\right\Vert ^{2},
\end{equation*}%
we deduce the required inequality (\ref{2.8}).

The sharpness of the constant may be proved as above in Theorem \ref{t2.1}
and we omit the details.
\end{proof}

\section{Another Reverse for the $\left( CBS\right) -$Inequality}

The following result also holds.

\begin{theorem}
\label{t3.1}Let $\left( X;\left( \cdot ,\cdot |\cdot \right) \right) $ be a
2-inner product space over $\mathbb{K}$ $\left( \mathbb{K}=\mathbb{C},%
\mathbb{R}\right) $ and $x,y,z\in X,$ $a,A\in \mathbb{K}$. If $A\neq -a$ and
either%
\begin{equation}
\func{Re}\left( Ay-x,x-ay|z\right) \geq 0  \label{3.1}
\end{equation}%
or, equivalently,%
\begin{equation}
\left\Vert x-\left. \frac{a+A}{2}y\right\vert z\right\Vert \leq \frac{1}{2}%
\left\vert A-a\right\vert \left\Vert y|z\right\Vert ,  \label{3.2}
\end{equation}%
holds, then we have the inequality%
\begin{align}
0& \leq \left\Vert x|z\right\Vert \left\Vert y|z\right\Vert -\func{Re}\left[ 
\func{sgn}\left( \frac{a+A}{2}\right) \left( x,y|z\right) \right]
\label{3.3} \\
& \leq \left\Vert x|z\right\Vert \left\Vert y|z\right\Vert -\left\vert
\left( x,y|z\right) \right\vert  \notag \\
& \leq \frac{1}{4}\frac{\left\vert A-a\right\vert ^{2}}{\left\vert
A+a\right\vert }\left\Vert y|z\right\Vert ^{2},  \notag
\end{align}%
where $\func{sgn}\left( \alpha \right) :=\frac{\alpha }{\left\vert \alpha
\right\vert },$ $\alpha \in \mathbb{C}\backslash \left\{ 0\right\} .$

The $\frac{1}{4}$ is best possible in the sense that it cannot be replaced
by a smaller constant.
\end{theorem}

\begin{proof}
We observe that the condition (\ref{3.2}) is equivalent with%
\begin{equation*}
\left\Vert x|z\right\Vert ^{2}-2\func{Re}\left[ \left( \frac{a+A}{2}\right)
\left( x,y|z\right) \right] +\left\vert \frac{a+A}{2}\right\vert
^{2}\left\Vert y|z\right\Vert ^{2}\leq \frac{1}{4}\left\vert A-a\right\vert
^{2}\left\Vert y|z\right\Vert ^{2}
\end{equation*}%
giving%
\begin{align}
\left\Vert x|z\right\Vert ^{2}+\left\vert \frac{a+A}{2}\right\vert
^{2}\left\Vert y|z\right\Vert ^{2}& \leq \frac{1}{4}\left\vert
A-a\right\vert ^{2}\left\Vert y|z\right\Vert ^{2}+2\func{Re}\left[ \left( 
\frac{a+A}{2}\right) \left( x,y|z\right) \right]  \label{3.4} \\
& \leq \frac{1}{4}\left\vert A-a\right\vert ^{2}\left\Vert y|z\right\Vert
^{2}+2\left\vert \frac{a+A}{2}\right\vert \left\vert \left( x,y|z\right)
\right\vert .  \notag
\end{align}%
By the elementary inequality%
\begin{equation*}
\alpha ^{2}+\beta ^{2}\geq 2\alpha \beta ,\ \ \ \alpha ,\beta \geq 0,
\end{equation*}%
we have%
\begin{equation}
2\left\vert \frac{a+A}{2}\right\vert \left\Vert x|z\right\Vert \left\Vert
y|z\right\Vert \leq \left\Vert x|z\right\Vert ^{2}+\left\vert \frac{a+A}{2}%
\right\vert ^{2}\left\Vert y|z\right\Vert ^{2}.  \label{3.5}
\end{equation}%
By making use of (\ref{3.4}) and (\ref{3.5}), we deduce%
\begin{align*}
0& \leq \left\vert \frac{a+A}{2}\right\vert \left\Vert x|z\right\Vert
\left\Vert y|z\right\Vert -\func{Re}\left[ \left( \frac{a+A}{2}\right)
\left( x,y|z\right) \right] \\
& \leq \left\vert \frac{a+A}{2}\right\vert \left[ \left\Vert x|z\right\Vert
\left\Vert y|z\right\Vert -\left\vert \left( x,y|z\right) \right\vert \right]
\\
& \leq \frac{1}{8}\left\vert A-a\right\vert ^{2}\left\Vert y|z\right\Vert
^{2},
\end{align*}%
which is clearly equivalent to the desired inequality (\ref{3.3}).

To prove the sharpness of the constant $\frac{1}{4}$ in (\ref{3.3}), let us
assume that there is a constant $D>0$ such that%
\begin{equation}
\left\Vert x|z\right\Vert \left\Vert y|z\right\Vert -\left\vert \left(
x,y|z\right) \right\vert \leq D\cdot \frac{\left\vert A-a\right\vert ^{2}}{%
\left\vert A+a\right\vert }\left\Vert y|z\right\Vert ^{2},  \label{3.6}
\end{equation}%
provided $x,y,z$ and $a,A$ satisfy the hypotheses of the theorem.

Assume now, $x,y,z,e\in X$ are such that $\left\Vert y,z\right\Vert =1,$ $%
\left\Vert e,z\right\Vert =1$ and $\left( e,y|z\right) =0.$ For $a,A\in 
\mathbb{K}$ with $a\neq -A,$ define%
\begin{equation*}
x=\frac{a+A}{2}y+\frac{A-a}{2}e.
\end{equation*}%
Then%
\begin{equation*}
\left\Vert x-\frac{a+A}{2},y|z\right\Vert =\frac{1}{2}\left\vert
A-a\right\vert ,
\end{equation*}%
and thus the condition (\ref{3.2}) is satisfied with equality.

Observe that, with the above choices for $x,y,z$ and $e$ we have%
\begin{align*}
\left\Vert x|z\right\Vert & =\sqrt{\frac{\left\vert A+a\right\vert ^{2}}{4}+%
\frac{\left\vert A-a\right\vert ^{2}}{4}}=\sqrt{\frac{\left\vert
A\right\vert ^{2}+\left\vert a\right\vert ^{2}}{2}}, \\
\left\vert \left( x,y|z\right) \right\vert & =\left\vert \frac{a+A}{2}%
\right\vert ,
\end{align*}%
and thus, from (\ref{3.6}), we deduce the inequality%
\begin{equation}
\sqrt{\frac{\left\vert A\right\vert ^{2}+\left\vert a\right\vert ^{2}}{2}}%
-\left\vert \frac{a+A}{2}\right\vert \leq D\cdot \frac{\left\vert
A-a\right\vert ^{2}}{\left\vert A+a\right\vert }  \label{3.7}
\end{equation}%
for $a,A\in \mathbb{C}$, $a\neq -A.$

For $\varepsilon \in \left( 0,1\right) ,$ consider $A=1+\sqrt{\varepsilon },$
$a=1-\sqrt{\varepsilon }.$ Then $a\neq -A$ and by (\ref{3.9}) we deduce%
\begin{equation*}
\sqrt{1+\varepsilon }-1\leq 2D\varepsilon ,
\end{equation*}%
giving by multiplication by $\sqrt{1+\varepsilon }+1>0$ that%
\begin{equation*}
\varepsilon \leq 2\varepsilon \left( \sqrt{1+\varepsilon }+1\right) D.
\end{equation*}%
Since $\varepsilon \in \left( 0,1\right) ,$ we may divide by $\varepsilon $
and thus we get%
\begin{equation}
D\geq \frac{1}{2\left( \sqrt{1+\varepsilon }+1\right) },\ \ \ \varepsilon
\in \left( 0,1\right) .  \label{3.8}
\end{equation}%
Letting $\varepsilon \rightarrow 0+$ in (\ref{3.8}), we obtain $D\geq \frac{1%
}{4},$ and the sharpness of the constant is proved.
\end{proof}

When the constants $A,a$ are real, we can point out the following reverse of
the triangle inequality.

\begin{corollary}
\label{c3.2}Let $\left( X;\left( \cdot ,\cdot |\cdot \right) \right) $ be a
2-inner product space over $\mathbb{K}$, $x,y,z\in X,$ and $m,M\in \left(
0,\infty \right) $ with $M>m.$ If either%
\begin{equation}
\func{Re}\left( My-x,x-my|z\right) \geq 0  \label{3.9}
\end{equation}%
or, equivalently,%
\begin{equation}
\left\Vert x-\left. \frac{m+M}{2}y\right\vert z\right\Vert \leq \frac{1}{2}%
\left( M-m\right) \left\Vert y|z\right\Vert  \label{3.10}
\end{equation}%
holds, then we have the inequality%
\begin{equation}
0\leq \left\Vert x|z\right\Vert +\left\Vert y|z\right\Vert -\left\Vert
x+y|z\right\Vert \leq \frac{1}{2}\cdot \frac{\left( M-m\right) }{\sqrt{M+m}}%
\left\Vert y|z\right\Vert .  \label{3.11}
\end{equation}
\end{corollary}

\begin{proof}
A simple computation shows that%
\begin{equation*}
\left( \left\Vert x|z\right\Vert +\left\Vert y|z\right\Vert \right)
^{2}-\left\Vert x+y|z\right\Vert ^{2}=2\left( \left\Vert x|z\right\Vert
\left\Vert y|z\right\Vert -\func{Re}\left( x,y|z\right) \right) .
\end{equation*}%
Using the inequality (\ref{3.3}), we may state that%
\begin{equation}
\left( \left\Vert x|z\right\Vert +\left\Vert y|z\right\Vert \right) ^{2}\leq
\left\Vert x+y|z\right\Vert ^{2}+\frac{1}{4}\frac{\left( M-m\right) ^{2}}{%
\left( M+m\right) }\left\Vert y|z\right\Vert ^{2}.  \label{3.12}
\end{equation}%
Taking the square root of (\ref{3.12}), we get%
\begin{align*}
\left\Vert x|z\right\Vert +\left\Vert y|z\right\Vert & \leq \sqrt{\left\Vert
x+y|z\right\Vert ^{2}+\frac{1}{4}\frac{\left( M-m\right) ^{2}}{\left(
M+m\right) }\left\Vert y|z\right\Vert ^{2}} \\
& \leq \left\Vert x+y|z\right\Vert +\frac{1}{2}\cdot \frac{\left( M-m\right) 
}{\sqrt{M+m}}\left\Vert y|z\right\Vert 
\end{align*}%
and the inequality (\ref{3.11}) is proved.
\end{proof}

\begin{remark}
\label{r3.3}Firstly, let us observe that from the inequality (\ref{1.15}) in
the Introduction, we may state the following additive reverse of the $\left(
CBS\right) -$inequality%
\begin{equation}
0\leq \left\Vert x|z\right\Vert \left\Vert y|z\right\Vert -\left\vert \left(
x,y|z\right) \right\vert \leq \frac{1}{2}\cdot \frac{\left\vert
A+a\right\vert -2\left[ \func{Re}\left( \bar{a}A\right) \right] ^{\frac{1}{2}%
}}{\left[ \func{Re}\left( \bar{a}A\right) \right] ^{\frac{1}{2}}}\left\vert
\left( x,y|z\right) \right\vert ,  \label{3.13}
\end{equation}%
provided $x,y,z\in X,$ $a,A\in \mathbb{K}$ with $\func{Re}\left( A\bar{a}%
\right) >0$ and either the condition (\ref{2.1}) or, equivalently (\ref{2.2}%
), is valid.

If $M>m>0$ and either (\ref{3.9}) or, equivalently, (\ref{3.10}) holds, then
from (\ref{3.13}) we may state the following simpler form%
\begin{equation}
0\leq \left\Vert x|z\right\Vert \left\Vert y|z\right\Vert -\left\vert \left(
x,y|z\right) \right\vert \leq \frac{1}{2}\cdot \frac{\left( \sqrt{M}-\sqrt{m}%
\right) ^{2}}{\sqrt{Mm}}\left\vert \left( x,y|z\right) \right\vert .
\label{3.14}
\end{equation}%
If, for the same $M,m$ we write the inequality (\ref{3.3}), then we have
another bound, namely:%
\begin{equation}
0\leq \left\Vert x|z\right\Vert \left\Vert y|z\right\Vert -\left\vert \left(
x,y|z\right) \right\vert \leq \frac{1}{4}\cdot \frac{\left( M-m\right) ^{2}}{%
\left( M+m\right) }\left\Vert y|z\right\Vert ^{2},  \label{3.15}
\end{equation}%
provided (\ref{3.9}), or equivalently, (\ref{3.10}) holds.
\end{remark}

\section{Integral Inequalities}

Let $\left( \Omega ,\Sigma ,\mu \right) $ be a measure space consisting of a
set $\Omega ,$ a $\sigma -$algebra $\Sigma $ of parts of $\Omega $ and a
countably additive and positive measure on $\Sigma $ with values in $\mathbb{%
R}\cup \left\{ \infty \right\} .$

Denote by $L_{\rho }^{2}\left( \Omega \right) ,$ the Hilbert space of all
real-valued functions $f$ defined on $\Omega $ that are $2-\rho -$integrable
on $\Omega .$ That is,%
\begin{equation*}
\int_{\Omega }\rho \left( t\right) \left\vert f\left( s\right) \right\vert
^{2}d\mu \left( s\right) <\infty ,
\end{equation*}%
where $\rho :\Omega \rightarrow \left( 0,\infty \right) $ is a measurable
function on $\Omega .$

If we denote by%
\begin{equation*}
\left\vert 
\begin{array}{lll}
a &  & b \\ 
&  &  \\ 
c &  & d%
\end{array}%
\right\vert ,\ \ \ \ \ \ \ \ a,b,c,d\in \mathbb{R}
\end{equation*}
the determinant associated with the matrix%
\begin{equation*}
\left( 
\begin{array}{lll}
a &  & b \\ 
&  &  \\ 
c &  & d%
\end{array}%
\right) ,\ \ \ \ \ \ \ \ a,b,c,d\in \mathbb{R};
\end{equation*}%
then we can introduce on $L_{\rho }^{2}\left( \Omega \right) $ the following
2-inner product%
\begin{multline}
\left( f,g|h\right) _{\rho }:=\frac{1}{2}\int_{\Omega }\int_{\Omega }\rho
\left( x\right) \rho \left( y\right) \left\vert 
\begin{array}{lll}
f\left( x\right)  &  & f\left( y\right)  \\ 
&  &  \\ 
h\left( x\right)  &  & h\left( y\right) 
\end{array}%
\right\vert   \label{4.1} \\
\times \left\vert 
\begin{array}{lll}
g\left( x\right)  &  & g\left( y\right)  \\ 
&  &  \\ 
h\left( x\right)  &  & h\left( y\right) 
\end{array}%
\right\vert d\mu \left( x\right) d\mu \left( y\right) ,
\end{multline}%
generating the 2-norm%
\begin{equation}
\left\Vert f|h\right\Vert _{\rho }=\left( \frac{1}{2}\int_{\Omega
}\int_{\Omega }\rho \left( x\right) \rho \left( y\right) \left\vert 
\begin{array}{lll}
f\left( x\right)  &  & f\left( y\right)  \\ 
&  &  \\ 
h\left( x\right)  &  & h\left( y\right) 
\end{array}%
\right\vert ^{2}d\mu \left( x\right) d\mu \left( y\right) \right) ^{\frac{1}{%
2}}.  \label{4.2}
\end{equation}%
A simple computation with integrals shows that%
\begin{equation*}
\left( f,g|h\right) _{\rho }=\left\vert 
\begin{array}{lll}
\dint_{\Omega }\rho \left( x\right) f\left( x\right) g\left( x\right) d\mu
\left( x\right)  &  & \dint_{\Omega }\rho \left( x\right) f\left( x\right)
h\left( x\right) d\mu \left( x\right)  \\ 
&  &  \\ 
\dint_{\Omega }\rho \left( x\right) g\left( x\right) h\left( x\right) d\mu
\left( x\right)  &  & \dint_{\Omega }\rho \left( x\right) h^{2}\left(
x\right) d\mu \left( x\right) 
\end{array}%
\right\vert 
\end{equation*}%
and%
\begin{equation*}
\left\Vert f|h\right\Vert _{\rho }=\left\vert 
\begin{array}{lll}
\dint_{\Omega }\rho \left( x\right) f^{2}\left( x\right) d\mu \left(
x\right)  &  & \dint_{\Omega }\rho \left( x\right) f\left( x\right) h\left(
x\right) d\mu \left( x\right)  \\ 
&  &  \\ 
\dint_{\Omega }\rho \left( x\right) f\left( x\right) h\left( x\right) d\mu
\left( x\right)  &  & \dint_{\Omega }\rho \left( x\right) h^{2}\left(
x\right) d\mu \left( x\right) 
\end{array}%
\right\vert ^{\frac{1}{2}}.
\end{equation*}%
We recall that the pair of functions $\left( q,p\right) \in L_{\rho
}^{2}\left( \Omega \right) \times L_{\rho }^{2}\left( \Omega \right) $ is
said to be \textit{synchronous} if%
\begin{equation*}
\left( q\left( x\right) -q\left( y\right) \right) \left( p\left( x\right)
-p\left( y\right) \right) \geq 0
\end{equation*}%
for a.e. $x,y\in \Omega .$

Now, suppose that $h\in L_{\rho }^{2}\left( \Omega \right) $ is such that $%
h\left( x\right) \neq 0$ for a.e. $x\in \Omega .$ Then by (\ref{4.1}) we
have the obvious identit, 
\begin{multline}
\left( f,g|h\right) _{\rho }=\frac{1}{2}\int_{\Omega }\int_{\Omega }\rho
\left( x\right) \rho \left( y\right) h^{2}\left( x\right) h^{2}\left(
y\right)   \label{4.3} \\
\times \left( \frac{f\left( x\right) }{h\left( x\right) }-\frac{f\left(
y\right) }{h\left( y\right) }\right) \left( \frac{g\left( x\right) }{h\left(
x\right) }-\frac{g\left( y\right) }{h\left( y\right) }\right) d\mu \left(
x\right) d\mu \left( y\right) 
\end{multline}%
and thus, a \textit{sufficient condition}\ for the inequality%
\begin{equation}
\left( f,g|h\right) _{\rho }\geq 0  \label{4.4}
\end{equation}%
to hold, is that the pair of functions $\left( \frac{f}{h},\frac{g}{h}%
\right) $ be synchronous. This condition is not necessary.

If $\Omega =\left[ a,b\right] \subset \mathbb{R}$ $\left( a<b\right) $ and $%
\mu $ is the Lebesgue measure, then a sufficient condition for the functions 
$\left( \frac{f\left( x\right) }{h\left( x\right) },\frac{g\left( x\right) }{%
h\left( x\right) }\right) ,$ $x\in \left[ a,b\right] $ to be synchronous is
that they are monotonic in the same sense, i.e. $\frac{f}{h}$ and $\frac{f}{g%
}$ are both increasing or decreasing on $\left[ a,b\right] .$ Obviously,
this condition is not necessary.

We are able now to state some integral inequalities that can be derived
using the general framework presented above.

\begin{proposition}
\label{p4.1}Let $M>m>0$ and $f,g,h\in L_{\rho }^{2}\left( \Omega \right) ,$ $%
h\neq 0,$ such that the functions%
\begin{equation}
M\cdot \frac{g}{h}-\frac{f}{h},\ \ \ \ \frac{f}{h}-m\cdot \frac{g}{h}
\label{4.5}
\end{equation}%
are synchronous on $\Omega .$ Then we have the inequalities%
\begin{eqnarray}
0 &\leq &\left\vert 
\begin{array}{lll}
\dint_{\Omega }\rho f^{2} &  & \dint_{\Omega }\rho fh \\ 
&  &  \\ 
\dint_{\Omega }\rho fh &  & \dint_{\Omega }\rho h^{2}%
\end{array}%
\right\vert \cdot \left\vert 
\begin{array}{lll}
\dint_{\Omega }\rho g^{2} &  & \dint_{\Omega }\rho gh \\ 
&  &  \\ 
\dint_{\Omega }\rho gh &  & \dint_{\Omega }\rho h^{2}%
\end{array}%
\right\vert   \label{4.6} \\
&&-\left\vert 
\begin{array}{lll}
\dint_{\Omega }\rho fg &  & \dint_{\Omega }\rho fh \\ 
&  &  \\ 
\dint_{\Omega }\rho gh &  & \dint_{\Omega }\rho h^{2}%
\end{array}%
\right\vert ^{2}  \notag \\
&\leq &\frac{1}{4}\left( M-m\right) ^{2}\left\vert 
\begin{array}{lll}
\dint_{\Omega }\rho g^{2} &  & \dint_{\Omega }\rho gh \\ 
&  &  \\ 
\dint_{\Omega }\rho gh &  & \dint_{\Omega }\rho h^{2}%
\end{array}%
\right\vert ^{2}  \notag \\
&&-\left\vert \frac{m+M}{2}\left\vert 
\begin{array}{lll}
\dint_{\Omega }\rho g^{2} &  & \dint_{\Omega }\rho gh \\ 
&  &  \\ 
\dint_{\Omega }\rho gh &  & \dint_{\Omega }\rho h^{2}%
\end{array}%
\right\vert -\left\vert 
\begin{array}{lll}
\dint_{\Omega }\rho fg &  & \dint_{\Omega }\rho fh \\ 
&  &  \\ 
\dint_{\Omega }\rho gh &  & \dint_{\Omega }\rho h^{2}%
\end{array}%
\right\vert \right\vert   \notag \\
&&\left( \leq \frac{1}{4}\left( M-m\right) ^{2}\left\vert 
\begin{array}{lll}
\dint_{\Omega }\rho g^{2} &  & \dint_{\Omega }\rho gh \\ 
&  &  \\ 
\dint_{\Omega }\rho gh &  & \dint_{\Omega }\rho h^{2}%
\end{array}%
\right\vert ^{2}\right) .  \notag
\end{eqnarray}
\end{proposition}

The proof is obvious by Theorem \ref{t2.1} and we omit the details.

The following counterpart of the $\left( CBS\right) -$inequality for
determinants also holds.

\begin{proposition}
\label{p4.2}With the assumptions of Proposition \ref{p4.1}, we have the
inequality%
\begin{eqnarray}
0 &\leq &\left\vert 
\begin{array}{lll}
\dint_{\Omega }\rho f^{2} &  & \dint_{\Omega }\rho fh \\ 
&  &  \\ 
\dint_{\Omega }\rho fh &  & \dint_{\Omega }\rho h^{2}%
\end{array}%
\right\vert \cdot \left\vert 
\begin{array}{lll}
\dint_{\Omega }\rho g^{2} &  & \dint_{\Omega }\rho gh \\ 
&  &  \\ 
\dint_{\Omega }\rho gh &  & \dint_{\Omega }\rho h^{2}%
\end{array}%
\right\vert  \label{4.7} \\
&&-\left\vert 
\begin{array}{lll}
\dint_{\Omega }\rho fg &  & \dint_{\Omega }\rho fh \\ 
&  &  \\ 
\dint_{\Omega }\rho gh &  & \dint_{\Omega }\rho h^{2}%
\end{array}%
\right\vert ^{2}  \notag
\end{eqnarray}%
\begin{eqnarray}
&\leq &\left( \frac{1}{4}\left( M-m\right) ^{2}\left\vert 
\begin{array}{lll}
\dint_{\Omega }\rho g^{2} &  & \dint_{\Omega }\rho gh \\ 
&  &  \\ 
\dint_{\Omega }\rho gh &  & \dint_{\Omega }\rho h^{2}%
\end{array}%
\right\vert \right.  \notag \\
&&-\left. \left\vert 
\begin{array}{lll}
\dint_{\Omega }\left( Mg-f\right) \left( f-mg\right) &  & \dint_{\Omega
}\rho \left( Mg-f\right) h \\ 
&  &  \\ 
\dint_{\Omega }\rho \left( f-mg\right) h &  & \dint_{\Omega }\rho h^{2}%
\end{array}%
\right\vert \right)  \notag \\
&&\times \left\vert 
\begin{array}{lll}
\dint_{\Omega }\rho g^{2} &  & \dint_{\Omega }\rho gh \\ 
&  &  \\ 
\dint_{\Omega }\rho gh &  & \dint_{\Omega }\rho h^{2}%
\end{array}%
\right\vert  \notag \\
&&\left( \leq \frac{1}{4}\left( M-m\right) ^{2}\left\vert 
\begin{array}{lll}
\dint_{\Omega }\rho g^{2} &  & \dint_{\Omega }\rho gh \\ 
&  &  \\ 
\dint_{\Omega }\rho gh &  & \dint_{\Omega }\rho h^{2}%
\end{array}%
\right\vert ^{2}\right) .  \notag
\end{eqnarray}
\end{proposition}

The proof follows by Theorem \ref{t2.2} applied for the 2-inner product
defined in (\ref{4.3}).

A different reverse of the $\left( CBS\right) -$inequality for determinants
is incorporated in the following proposition.

\begin{proposition}
\label{p4.3}With the assumptions of Proposition \ref{p4.1}, we have the
inequality%
\begin{eqnarray}
0 &\leq &\left\vert 
\begin{array}{lll}
\dint_{\Omega }\rho f^{2} &  & \dint_{\Omega }\rho fh \\ 
&  &  \\ 
\dint_{\Omega }\rho fh &  & \dint_{\Omega }\rho h^{2}%
\end{array}%
\right\vert ^{\frac{1}{2}}\cdot \left\vert 
\begin{array}{lll}
\dint_{\Omega }\rho g^{2} &  & \dint_{\Omega }\rho gh \\ 
&  &  \\ 
\dint_{\Omega }\rho gh &  & \dint_{\Omega }\rho h^{2}%
\end{array}%
\right\vert ^{\frac{1}{2}}  \label{4.8} \\
&&-\left\vert \det \left( 
\begin{array}{lll}
\dint_{\Omega }\rho fg &  & \dint_{\Omega }\rho fh \\ 
&  &  \\ 
\dint_{\Omega }\rho gh &  & \dint_{\Omega }\rho h^{2}%
\end{array}%
\right) \right\vert   \notag \\
&\leq &\frac{1}{4}\frac{\left( M-m\right) ^{2}}{M+m}\left\vert 
\begin{array}{lll}
\dint_{\Omega }\rho g^{2} &  & \dint_{\Omega }\rho gh \\ 
&  &  \\ 
\dint_{\Omega }\rho gh &  & \dint_{\Omega }\rho h^{2}%
\end{array}%
\right\vert ^{2}.  \notag
\end{eqnarray}%
The constant $\frac{1}{4}$ is best possible in (\ref{4.8}).
\end{proposition}

The proof follows from Theorem \ref{t3.1} applied for the 2-inner product
defined in (\ref{4.3}).

Finally, by the use of Corollary \ref{c3.2}, we may state the following
reverse of the triangle inequality for determinants.

\begin{proposition}
\label{p4.4}With the assumptions of Proposition \ref{p4.1}, we have the
inequality:%
\begin{eqnarray}
0 &\leq &\left\vert 
\begin{array}{lll}
\dint_{\Omega }\rho f^{2} &  & \dint_{\Omega }\rho fh \\ 
&  &  \\ 
\dint_{\Omega }\rho fh &  & \dint_{\Omega }\rho h^{2}%
\end{array}%
\right\vert ^{\frac{1}{2}}+\left\vert 
\begin{array}{lll}
\dint_{\Omega }\rho gh &  & \dint_{\Omega }\rho gh \\ 
&  &  \\ 
\dint_{\Omega }\rho gh &  & \dint_{\Omega }\rho h^{2}%
\end{array}%
\right\vert ^{\frac{1}{2}}  \label{4.9} \\
&&-\left\vert 
\begin{array}{lll}
\dint_{\Omega }\rho \left( f+g\right) ^{2} &  & \dint_{\Omega }\rho \left(
f+g\right) h \\ 
&  &  \\ 
\dint_{\Omega }\rho \left( f+g\right) h &  & \dint_{\Omega }\rho h^{2}%
\end{array}%
\right\vert ^{\frac{1}{2}}  \notag \\
&\leq &\frac{1}{2}\cdot \frac{M-m}{\sqrt{M+m}}\cdot \left\vert 
\begin{array}{lll}
\dint_{\Omega }\rho g^{2} &  & \dint_{\Omega }\rho gh \\ 
&  &  \\ 
\dint_{\Omega }\rho gh &  & \dint_{\Omega }\rho h^{2}%
\end{array}%
\right\vert ^{\frac{1}{2}}.  \notag
\end{eqnarray}
\end{proposition}

\textbf{Acknowledgement:}{\small \ S. S. Dragomir and Y. J. Cho greatly
acknowledge the financial support from the Brain Pool Program (2002) of the
Korean Federation of Science and Technology Societies. The research was
performed under the "Memorandum of Understanding" between Victoria
University and Gyeongsang National University.}

\end{document}